\documentclass[a4paper, 11pt]{article}

\usepackage[utf8]{inputenc} 
\usepackage[T1]{fontenc}    
\usepackage{hyperref}       
\usepackage{url}           
\usepackage{booktabs}      
\usepackage{amsfonts}     
\usepackage{nicefrac}       
\usepackage{microtype}     
\usepackage{lipsum}
\usepackage{tikz}
\usepackage[noend]{algpseudocode}
\usepackage{algorithmicx, algorithm}
\usepackage{xcolor}
\usepackage[misc]{ifsym}
\usepackage{graphicx}
\usepackage{epstopdf}
\usepackage{subfigure}
\usepackage{amsmath}
\usepackage{todonotes}
\usepackage{comment}
\usepackage{amsthm}
\newtheorem*{remark}{Remark}
\newcommand\LL[1]{{\color{black}#1}}

\title{Finite Difference Nets: A Deep Recurrent Framework for Solving Evolution PDEs}

 \author{
 Cheng Chang \thanks{Department of Mathematics and Institute of Mathematical Sciences, The Chinese University of Hong Kong (chengchang@link.cuhk.edu.hk)},~
 Liu Liu \thanks{Department of Mathematics, The Chinese University of Hong Kong (lliu@math.cuhk.edu.hk)}, ~Tieyong Zeng \thanks{Department of Mathematics, The Chinese University of Hong Kong (zeng@math.cuhk.edu.hk)}
 }

\begin{document}
\maketitle

\begin{abstract}

\LL{There has been an arising trend of adopting deep learning methods to study partial differential equations (PDEs).
In this paper, we introduce a deep recurrent framework for solving time-dependent PDEs without generating large scale data sets.
We provide a new perspective, that is, a different type of architecture
through exploring the possible connections between traditional numerical methods (such as finite difference schemes) and deep neural networks, particularly convolutional and fully-connected neural networks are employed. Our proposed approach will show its effectiveness and efficiency in solving PDE models with an integral form, in particular, we test on one-way wave equations and system of conservation laws.
}
\end{abstract}


\section{Introduction}

\LL{In recent years, deep learning, or machine learning in general, has drawn an incredibly arising attention and achieved significant developments in various research fields, such as image recognition \cite{Simonyan2014}, image synthesis \cite{Odena2017}, natural language processing \cite{Gardner2018} and speech recognition \cite{Amodei2015}. There has been quite a few work of adopting deep learning approaches to numerically solve partial differential equations (PDEs). 
For many types of PDEs, despite robust numerical methods have been developed, for numerical analyst there have always been big challenges in designing efficient traditional solvers.

The deep learning method, on the other hand, provides a new and meshfree approach for solving partial differential equations (PDEs) that can resolve many difficulties appeared in traditional numerical methods \cite{raissi2019physics}, such as the expensive computational cost for high-dimensional problems, the challenges dealing with complex boundary conditions, truncation of the computational domain, etc. Compared to conventional numerical methods that may require a solid understanding on numerical analysis, deep learning algorithms have the advantage of being intuitive and easy to be executed. We mention that there are other meshfree methods to solve PDEs, such as by using evolutionary algorithms to search for the optimized solution, in particular elliptic type equations were studied in \cite{EA2014}. Besides, there have been recent works on using deep learning methods to study other types of equations, for example \cite{DNN-LTE, AP-NN, DNN-LTE1, FP-NN}.}

The highlights of our contributions are as the following. First, we introduce a deep learning method to find solutions to the evolution equations, without involving large data sets. Second, 
compared to several existing literature that studied deep learning methods on this topic, we provide a different type of architecture through exploring the possible connection between finite difference methods and deep neural network approaches that particularly use convolutional and fully-connected layers. 
Thirdly, our proposed method is naturally conducive to solving equation with an integral form, such as the conservation laws. Furthermore, we observe in the numerical experiments that adopting the integral form seems to be more accurate than the differential form. 

\LL{Besides, we also devote some efforts to solve inverse problems. The goal is to explore the underlying differential equation models with some training data given. Note that in \cite{Qin2019}, the authors studied a numerical framework for approximating unknown governing equations using observation data and deep neural networks. They demonstrate that the residual network and its extension can be considered as one-step or multi-step methods that are exact in temporal integration. One main difference of our work is that no data (except for possible initial and boundary conditions) is required. There are many relevant work in this direction \cite{Inverse-design, B-PINN, ZGG20}.}

\LL{Another series of very successful and popular methods are based on the physics-informed neural network (PINN) framework, where broad applications arise in industry and engineering \cite{han2018solving, jin2020nsfnets, kissas2020machine, raissi2019physics, raissi2020hidden}. It leverages the benefits of auto-differentiation of the current software and the underlying physics in the models, incorporates it into the network by minimizing the loss function composed of the PDE's residuals and other constraints such as the initial and boundary conditions.
It is also worth mentioning that in \cite{RNN}, a recursive deep neural network approach has been developed to study dynamical systems, and retrieve the governing equation accurately, under weak requirements on the existing data, such as sampled with big time lags.}

\LL{This paper is organized by the following. In Section \ref{sec:Setup}, we introduce the problem setup and our deep recurrent framework for solving evolution equations. Our proposed method takes solution values at the current time step as inputs and produces solution values at the next time step, thus theoretically it mimics the evaluation of finite difference schemes. In Section \ref{sec:Num}, several numerical experiments for different PDEs will be conducted to demonstrate the performance of our deep learning approach. Finally, we conclude the paper in Section \ref{sec:Con}. }

\section{Our method}
\label{sec:Setup}

\subsection{Convolutional neural network}

\LL{{\bf Motivation.

}
There has been numerous study on using finite difference schemes that possess different properties to solve time-dependent PDEs \cite{LeVeque1992}. 
To summarize, we approximate the solution at current time step $u^{(n)}(x)$ by a linear combination of solutions obtained at previous time step and in some neighbourhood of $x$, which is denoted by $N(x)$. 
For example, for first-order in time scheme, one computes
\begin{equation}
\label{FD}
u^{(n)}(x)=\sum_{x_i\in N(x)} \alpha_i u^{(n-1)}(x_i). 
\end{equation}
Here, the neighbourhood $N(x)$ and coefficients $\alpha_i$ can be determined according to the finite difference scheme used. 
}

\LL{Deep learning, as one of the subsets of machine learning, uses deep learning algorithms to implicitly extract important information based on input data. It is designed to approximate a unknown target function by using the deep neural network, for example see an overview \cite{Goodfellow2016, DL-Review}.} Deep neural networks compose computations performed by many layers, within which each layer consists of a linear mapping (e.g. matrix multiplication in fully-connected layer), a (usually nonlinear) activation function, and a loss function whose gradient is computed then backpropagated to the first layer in order to update all the parameters in the neural network. Then this procedure will be repeated until the loss function converges to zero. 

\LL{Several types of neural networks formulate the main architectures of deep learning. Among them, a popular used and important architecture is the convolutional neural network (CNN), which 
we found it reminiscent of the above formulation \eqref{FD}.}
The computation of a typical convolutional layer with kernel function $K$ and input $I$ can be written as
\begin{equation}
\text{Conv}(x,I;K)=\sum_{y\in N(x)} I(x)K(y-x). 
\end{equation}
Here the kernel values can be thought of as the constant coefficients in the finite difference scheme. 
\LL{Nevertheless, we do not need to artificially design them, since appropriate values can be found automatically via the optimization algorithm adopted. }

\begin{figure}
  \centering
  \includegraphics[width=0.6\textwidth]{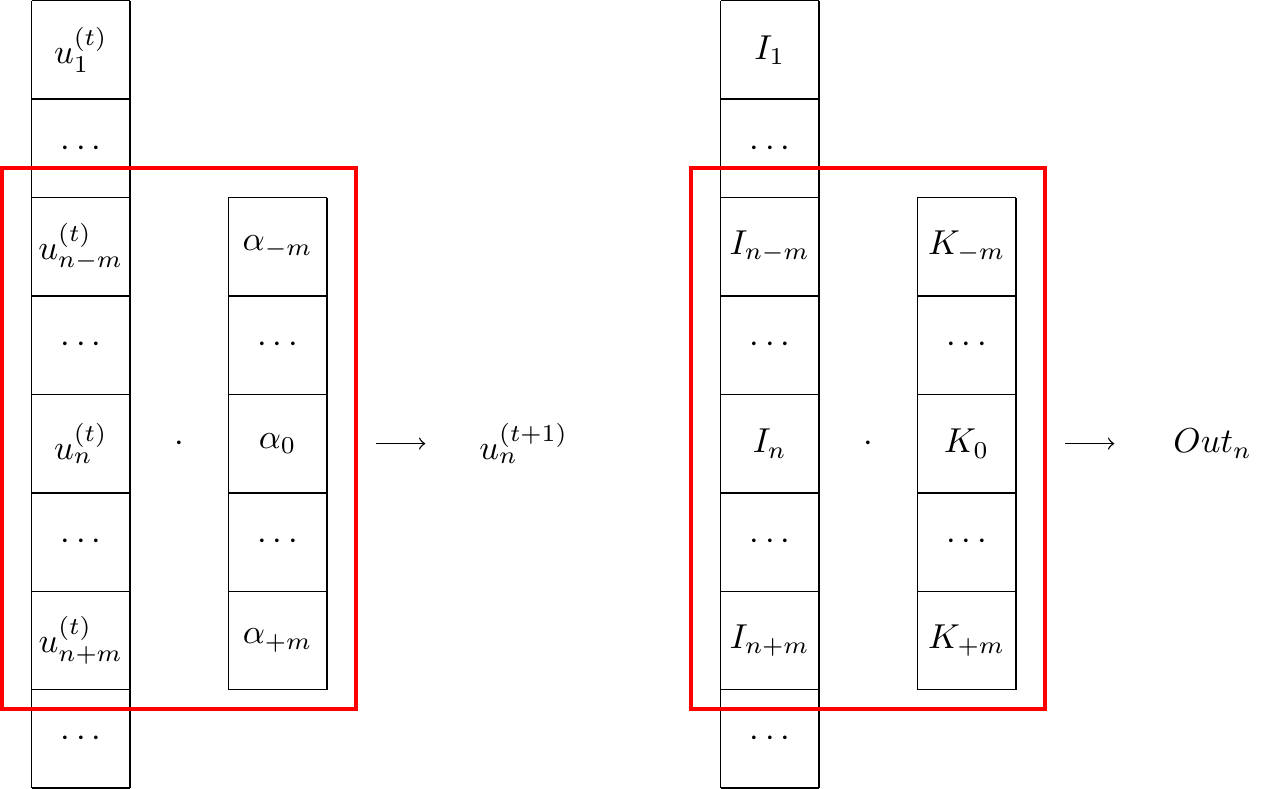}
  \caption{Similarity between finite-difference schemes and convolutional layers.}
  \label{fig:fig1}
\end{figure}

We remark that it is also reasonable to use fully-connected layers. The calculation of a convolutional layer can be reformulated as a multiplication of a matrix with special requirement on its structure, which indicates that a fully-connected layer is able to handle more general cases at an acceptable computational cost \cite{Dumoulin2018}. 

%--------------------------------
\subsection{Conservation laws}

\LL{In this paper, we consider hyperbolic systems of conservation laws \cite{LeVeque1992}. In particular, we will numerically study the time-dependent systems of PDEs in one space dimension that take the form: 
\begin{equation}
\label{PDE-model}
\frac{\partial}{\partial t}u(x,t) + \frac{\partial}{\partial x}F(u(x,t))=0, 
\end{equation}
where $u: \mathbb{R}\times \mathbb{R} \to \mathbb{R}^m$ is an $m$-dimensional vector of conserved quantities
or state variables, such as mass, momentum and energy in a fluid dynamics problem. 
The vector-valued function $F(u)$ is called the flux function for the system of conservation laws, and its form is known.} The initial condition for \eqref{PDE-model} is given by
\begin{equation}
u(x, 0)=u_0(x), 
\end{equation}
and the boundary condition is 
\begin{equation}
u(x, t)=u_B(x, t), \qquad \text{for } x\in \partial\Omega, 
\end{equation}
with functions $u_0$ and $u_B$ given. 

%--------------------------------
\subsection{Proposed Framework}

\subsubsection{A differential Form}

Based on the above observation, a fully-connected or convolutional layer can be treated as an one time step advancement in different finite difference schemes. 
\LL{We divide the spatial domain $\Omega$ into a uniform mesh.}
The input of the neural network is a tensor $u^{(n)}$ storing the function values of the current time step at all the spatial grid points. The output is approximated function values at the next time step $u^{(n+1)}$. 
For the model \eqref{PDE-model}, we design the loss function accordingly: 
\begin{equation}
\label{Loss}
\text{Loss}(\mathcal N(u^{(n)}))=\int_{\Omega}\big |\frac{\partial u^{(n)}}{\partial t}- \frac{\partial}{\partial x} F(\mathcal{N}(u^{(n)})\big |^2\, dx+\int_{\partial\Omega}\big |u^{(n)}(x,t)-u_B(x,t)\big |^2\, dS, 
\end{equation}
\LL{where the operator $\mathcal N$ is used to represent the processing of the neural network, of which output denoted by $\mathcal N(u^{(n)})$.}

Notice that the loss function is designed similarly to that in the PINN framework \cite{Raissi2019}, except that the initial condition needs not to be incorporated here. 
The time derivative can be approximated by the difference of input and output of the network divided by time step size $\Delta t$: 
\begin{equation}
\frac{\partial u^{(n)}}{\partial t}\approx \frac{\mathcal{N}(u^{(n)})-u^{(n)}}{\Delta t}. 
\end{equation}
Similarly, the derivative of flux $F$ with respect to the spatial variable can be approximated by using its
neighbouring points: 
\begin{equation}
\frac{\partial}{\partial x} F(\mathcal{N}(u^{(n)})\big |_{x=x_0}\approx 
\frac{F(\mathcal{N}(u^{(n)}(x_0))) - F(\mathcal{N}(u^{(n)}(x_0-\Delta x)))}{\Delta x}.
\end{equation}
\LL{If there are higher-order derivatives involved in the model equation, the same procedure can be employed to calculate them. We also remark that in the definition of loss function \eqref{Loss}, the square of $l_2$-norm is adopted. Instead, if one uses simply the $l_2$-norm itself, there will be numerical instability due to the $l_2$-norm values shown in the denominator of the gradients \cite{github.com}.}
Figure \ref{fig:overview} below illustrates an overview of our proposed framework.

\begin{figure}
  \centering
  \includegraphics[width=0.45\textwidth]{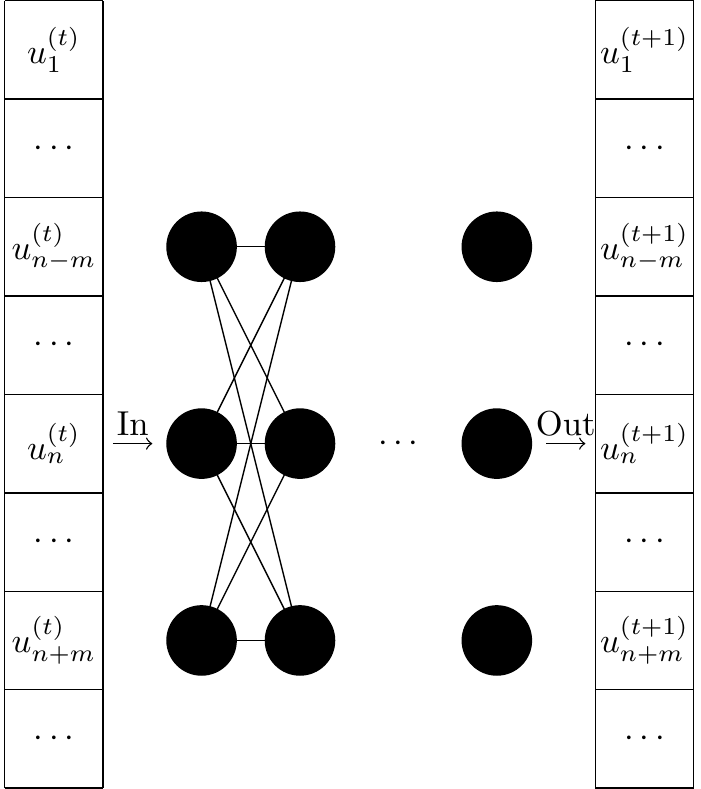}
  \caption{Overview of the proposed method.}
  \label{fig:overview}
\end{figure}

We now summarize the training procedure in {\bf Algorithm 1} as follows. 
\begin{algorithm}[!t]
\caption{Training Procedure}
\hspace*{0.02in}{\bf Input: } Let the time and space step size be $\Delta t$ and $\Delta x$. Denote the number of iterations per time step as $n$. \\
\hspace*{0.02in}{\bf Output: } Function values for $u$ at each time step.
\begin{algorithmic}[1]
\State Divide $\Omega$ into small grids by using $\Delta x$. 
\State Let $I$ be the tensor storing the initial values $u_0$ at all spatial grid points. 
\State $t\leftarrow 0$
\While{$t<T$}
	\State $t\leftarrow t+\Delta t$
	\State $i\leftarrow 0$
	\While{$i<n$}
		\State $I_{\text{new}}\leftarrow \mathcal{N}(I)$
		\State $\text{Loss}\leftarrow \text{Loss}(I_{\text{new}})$
		\State Apply gradients to the neural network.
	\EndWhile
	\State Return $I$ as the function value at time $t$. 
	\State $I\leftarrow I_{\text{new}}$. 
\EndWhile
\end{algorithmic}
\end{algorithm}

\bigskip

\subsubsection{An integral Form}
\LL{In the above framework, since the function values at all mesh points are readily available, it becomes straightforward to compute based on integration of the model \eqref{PDE-model} (in $t$ and $x$, on both sides). 
For a general function $g$, one can get its integration in $t$ and $x$ by applying the trapezoidal rule.}
\begin{equation*}
\int_{t_k}^{t_{k+1}}g(x, t)\,dt\approx\frac12(g(x, t_k)+g(x, t_{k+1}))\Delta t\,,
\end{equation*}
\begin{equation*}
\int_{x_k}^{x_{k+m}}g(x, t)\,dx\approx\frac12\sum_{i=k}^{k+m-1}(g(x_i, t)+g(x_{i+1}, t))\Delta x\,. 
\end{equation*}

It is worth mentioning that the integral form seems more accurate than the differential form presented earlier, from our numerical experiments. This is due to the fact that the numerical method we used to approximate the derivatives is only first-order accurate, while the trapezoidal rule adopted in the integration is of second-order accuracy. 
\LL{We will investigate how to improve the accuracy of neural network approximations in our future work.}
Other aspects of this method using the integral form, such as the training procedure and loss function, are the same as introduced in the differential form.

\begin{remark}
\LL{We make the remark that our proposed method not only works for solving the conservation laws but can also be applied to general high-dimensional, time-dependent PDEs. For simplicity, we illustrate our idea by studying two simple one-dimensional fluid models in Section \ref{sec:Num}. We defer to our future work on testing other types of higher dimensional PDEs.}
\end{remark}

%-------------------------------------------------------------------
\section{Numerical Examples}
\label{sec:Num}

\subsection{One-way Wave Equation}

We first test our proposed method on the one-way wave equation given by
\begin{equation}
\frac{\partial u}{\partial t}-\frac{\partial u}{\partial x}=0, \qquad (x,t)\in[0,1]\times [0,1],
\end{equation}
with initial condition $u_0(x)=0$ and Dirichlet boundary condition:
\begin{equation}
u_B(1, t)=
\left\{
\begin{array}{c}
\sin(\pi t)^2, \quad t\leq 1, \\
0, \qquad\qquad t>1. 
\end{array}
\right.
\end{equation}

\begin{figure}[!h]
  \centering
  \subfigure{
  \begin{minipage}[b]{.3\linewidth}
  \includegraphics[width=1\textwidth,height=0.9\textwidth]{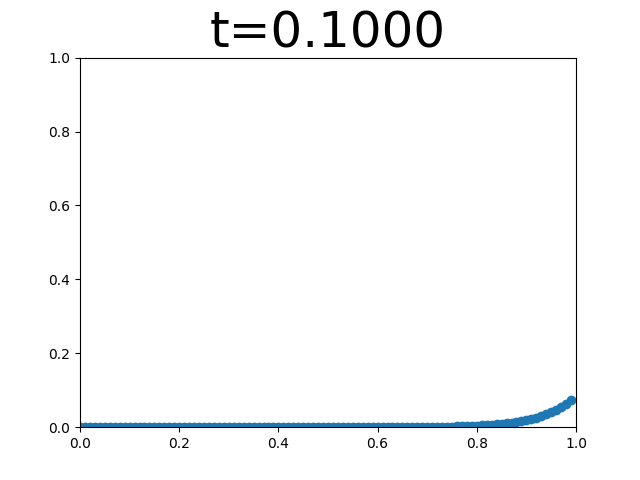}
  \end{minipage}
  }
  \subfigure{
  \begin{minipage}[b]{.3\linewidth}
  \includegraphics[width=1\textwidth,height=0.9\textwidth]{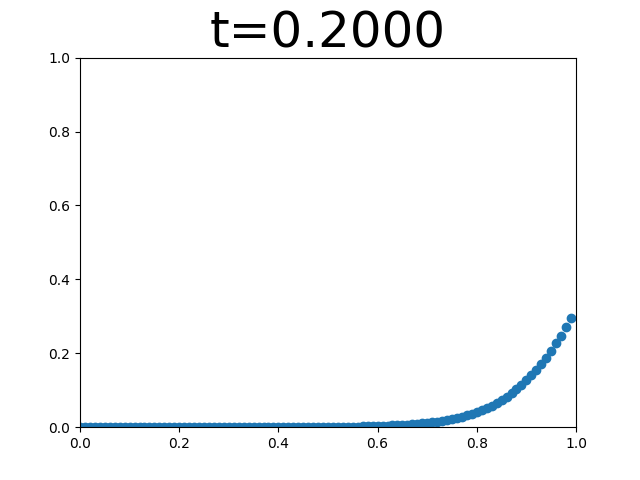}
  \end{minipage}
  }
  \subfigure{
  \begin{minipage}[b]{.3\linewidth}
  \includegraphics[width=1\textwidth,height=0.9\textwidth]{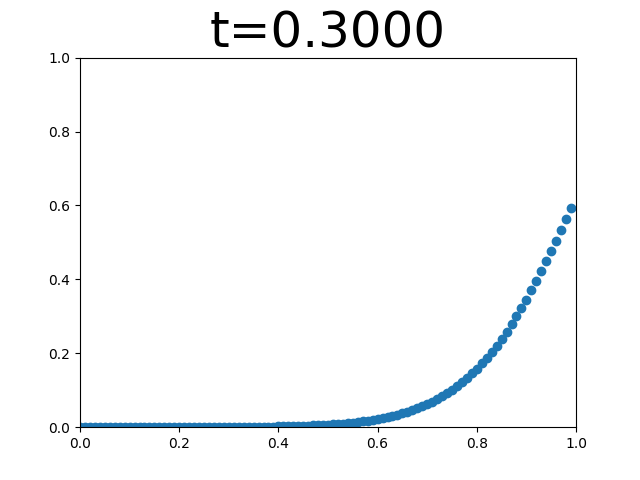}
  \end{minipage}
  }
  \subfigure{
  \begin{minipage}[b]{.3\linewidth}
  \includegraphics[width=1\textwidth,height=0.9\textwidth]{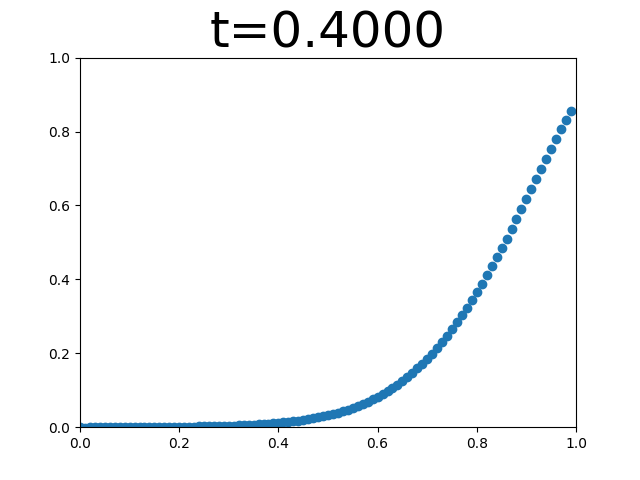}
  \end{minipage}
  }
  \subfigure{
  \begin{minipage}[b]{.3\linewidth}
  \includegraphics[width=1\textwidth,height=0.9\textwidth]{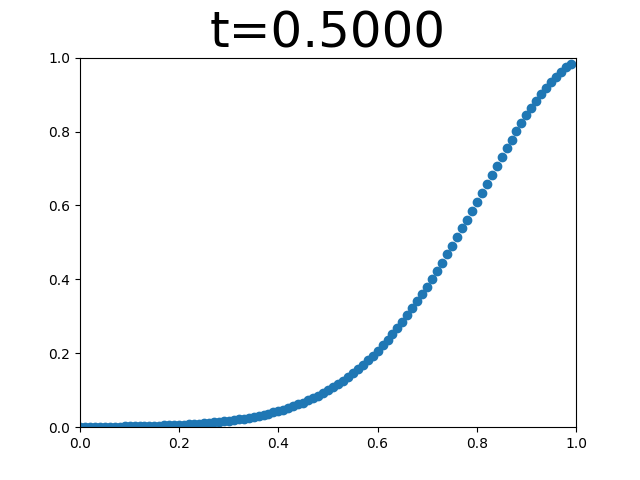}
  \end{minipage}
  }
  \subfigure{
  \begin{minipage}[b]{.3\linewidth}
  \includegraphics[width=1\textwidth,height=0.9\textwidth]{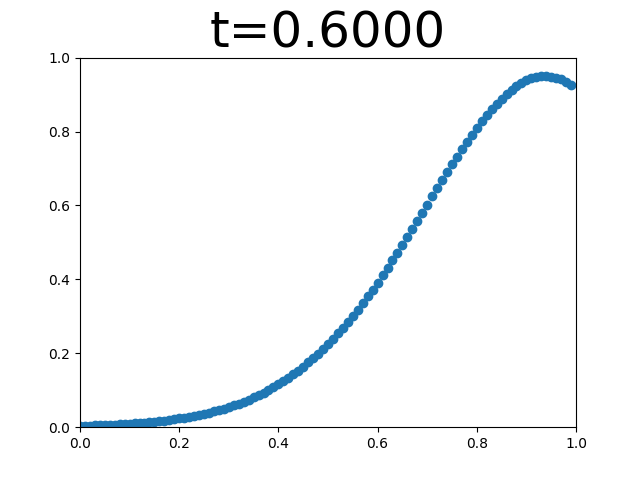}
  \end{minipage}
  }
  \subfigure{
  \begin{minipage}[b]{.3\linewidth}
  \includegraphics[width=1\textwidth,height=0.9\textwidth]{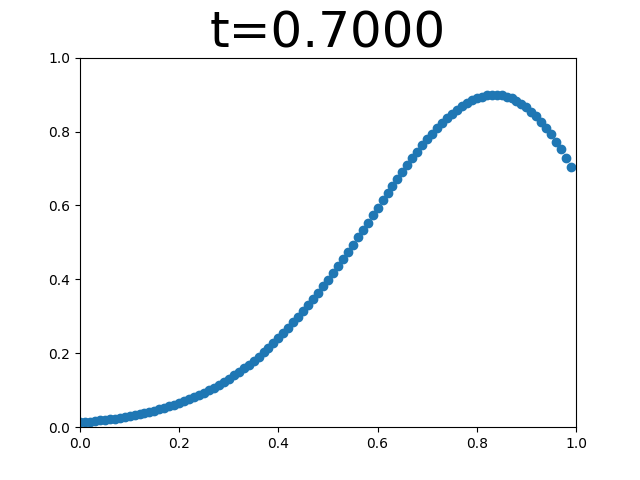}
  \end{minipage}
  }
  \subfigure{
  \begin{minipage}[b]{.3\linewidth}
  \includegraphics[width=1\textwidth,height=0.9\textwidth]{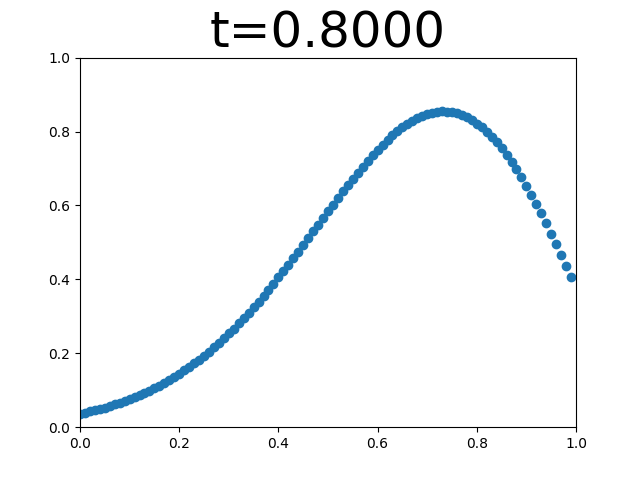}
  \end{minipage}
  }
  \subfigure{
  \begin{minipage}[b]{.3\linewidth}
  \includegraphics[width=1\textwidth,height=0.9\textwidth]{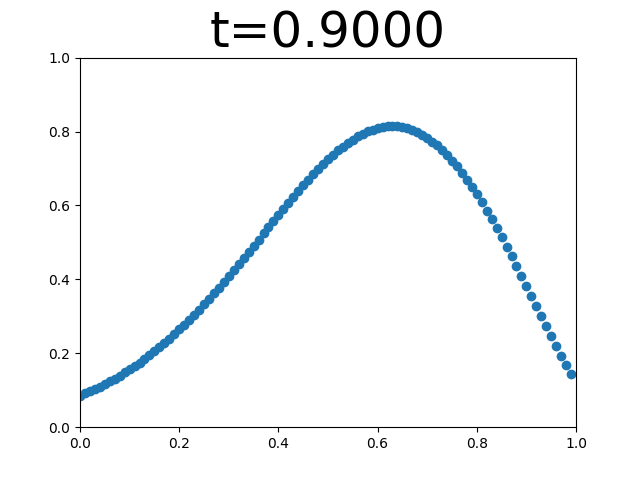}
  \end{minipage}
  }
  \label{fig:wave}
  \caption{Neural network trained solutions for the one-way wave equation. }
\end{figure}

The gap size between the sampling points is $0.01$. Our trained network contains $3$ fully-connected layers with $100$ hidden units per layer. Both non-output layers are activated using rectified linear units. The numerical results are shown in Figure \ref{fig:wave}. One can observe from solution approximations that a wave propagates from the right to the left. This phenomenon is consistent with the behaviour of the solution of the one-way wave equation.

\subsection{Euler Equations}

\LL{In fluid dynamics, the Euler equations are a set of quasilinear hyperbolic equations governing adiabatic and inviscid flow. For compressible Euler equations, please refer to the review articles \cite{Euler2,Euler1}. In one space dimension, the conservative form is written as}
\begin{equation}
\label{Euler}
\textcolor{black}{\frac{\partial U}{\partial t}+ \frac{\partial F}{\partial x}=0, }
\end{equation}
where the conserved quantity and associated flux are
\begin{equation}
U=
\left(
\begin{array}{c}
\rho \\
\rho u \\
E
\end{array}
\right), \qquad F=
\left(
\begin{array}{c}
\rho u \\
\rho u^2 + p \\
(E+p)u
\end{array}
\right). 
\end{equation}
\LL{Here
\begin{equation}
    \label{PP}
    p = \rho T
\end{equation}
is the pressure, and the relation between the energy $E$ and temperature $T$ is given by
\begin{equation}
\label{E-T}
E = \frac{1}{2}\rho u^2 + \rho T. \end{equation}}

In our numerical test, we consider the shock tube problem \cite{Sod1978} on the domain $\Omega=[0, 1]$, with initial conditions: 
\begin{equation}
\rho_0(x)=
\left\{
\begin{array}{rcl}
8, \quad x\in [0, 0.5) \\
1,\quad x\in [0.5, 1]
\end{array}
\right., \quad
u_0(x) = 0, \quad
p_0(x) = 
\left\{
\begin{array}{rcl}
8, \quad x\in [0, 0.5)\\
1, \quad x\in [0.5, 1]
\end{array}\right.. 
\end{equation}
Note that the corresponding initial conditions for $E$ can be obtained from \eqref{PP}--\eqref{E-T}. 
Dirichlet boundary conditions are assumed. 

\bigskip
\bigskip

\noindent$\bullet${\textbf{   Integral Form}}
\\[2pt]
\indent \LL{Applying integration in $t$ and $x$, the integral form of the Euler's equations is given by}
\begin{equation}
\int_{x_1}^{x_2}u(x, t_2)\,dx=\int_{x_1}^{x_2}u(x, t_1)\,dx+
\int_{t_1}^{t_2}F(u(x_1, t))\,dt-\int_{t_1}^{t_2}F(u(x_2, t))\,dt
\end{equation}
for any $x_1, x_2 \in \Omega$ and $t_1$, $t_2$. 

In our proposed framework, the time $t_1$ corresponds to the network input in one iteration, while the time $t_2$ is associated with the network output. The equation is enforced to be valid on every small subinterval, that being said, we let $x_1$ and $x_2$ be two neighbouring grid points.
The derivation from the true solution at the small grid $[x_k, x_{k+1}]$ is shown as
\begin{equation}
  \begin{aligned}
   l_k =  & \, \frac12 \left(\mathcal{N}(u^{(n)})\big |_{x=x_k} + \mathcal{N}(u^{(n)})\big |_{x=x_{k+1}}\right)
  \Delta x - \frac12\left(u^{(n)}\big |_{x=x_k}+u^{(n)}\big |_{x=x_{k+1}}\right)\Delta x \\[4pt] 
  & - \frac12 \left(F(u^{(n)})\big |_{x=x_k}+F(\mathcal{N}(u^{(n)}))\big |_{x=x_k}\right)\Delta t \\[4pt]
  & + \frac12 \left(F(u^{(n)})\big |_{x=x_{k+1}}+F(\mathcal{N}(u^{(n)}))\big |_{x=x_{k+1}}\right)\Delta t. \end{aligned}
\end{equation}
The overall loss function regarding the integral of the model equation is defined by
\begin{equation}
l=\sum_k l^2_k+l_{BC}, 
\end{equation}
where $l_{\text{BC}}$ is the same as the second term in \eqref{Loss}. 
\bigskip
\bigskip 

\noindent $\bullet${\textbf{  Initialization of parameters}}
\\[2pt]
\indent Although rarely used in deep learning practices, we find that initializing the parameters in our neural network with identity matrix is preferable. The reason may be due to that in the physics system we studied, the evolution of state variables is usually smooth in time, so the optimal parameters in the network are not expected to be far away from identity. On the other hand, initializing the parameters with identities seem can also accelerate the convergence of neural network and avoid getting trapped by local minima.  

\subsection{Numerical Results for the Euler system}

\LL{It is worth mentioning that in \cite{PINN-Flow}, the authors have employed the PINNs framework to approximate the Euler equations that model high-speed aerodynamic flows; both the forward and inverse problems in one-dimensional and two-dimensional domains were studied carefully. Here we numerically study the compressible Euler system \eqref{Euler} by using our proposed approach.}

In this test, the gap size between the sampling points is set $0.01$ as well. The network contains one hidden layer with $330$ hidden units. The input includes three physical quantities, so there will be a $3\times 100$ matrix passing a flatten layer to obtain a $300$-dimensional vector. The layer parameters are initialised to be an identity matrix. 

\begin{figure}[!h]
  \centering
  \subfigure{
  \begin{minipage}[b]{.4\linewidth}
  \includegraphics[width=1.1\textwidth,height=0.9\textwidth]{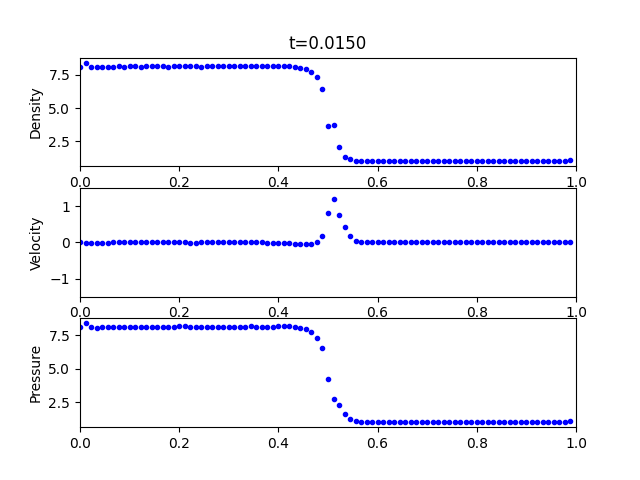}
  \end{minipage}
  }
  \subfigure{
  \begin{minipage}[b]{.4\linewidth}
  \includegraphics[width=1.1\textwidth,height=0.9\textwidth]{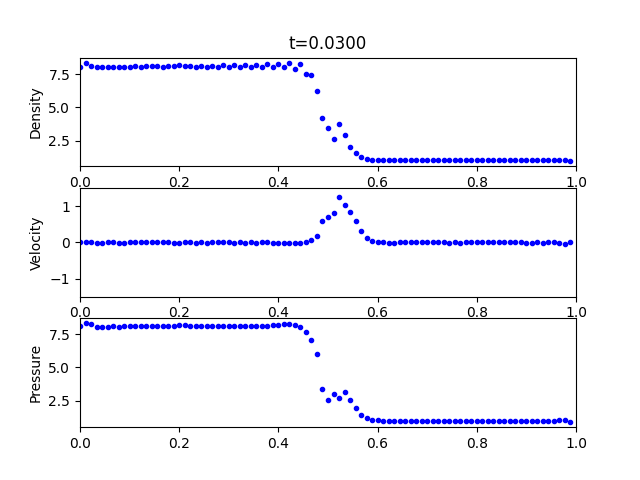}
  \end{minipage}
  }
  \subfigure{
  \begin{minipage}[b]{.4\linewidth}
  \includegraphics[width=1.1\textwidth,height=0.9\textwidth]{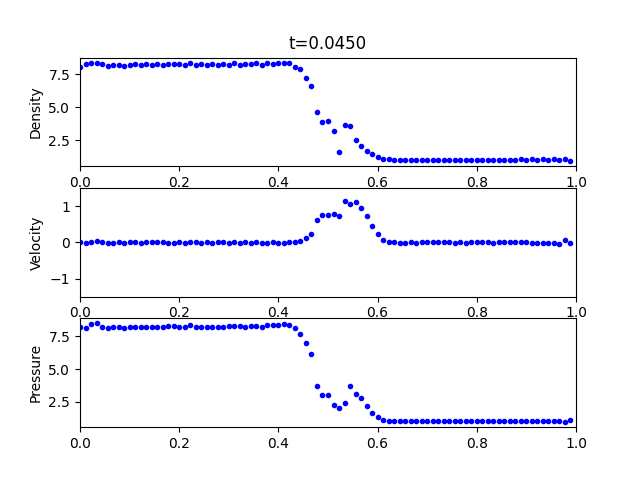}
  \end{minipage}
  }
  \subfigure{
  \begin{minipage}[b]{.4\linewidth}
  \includegraphics[width=1.1\textwidth,height=0.9\textwidth]{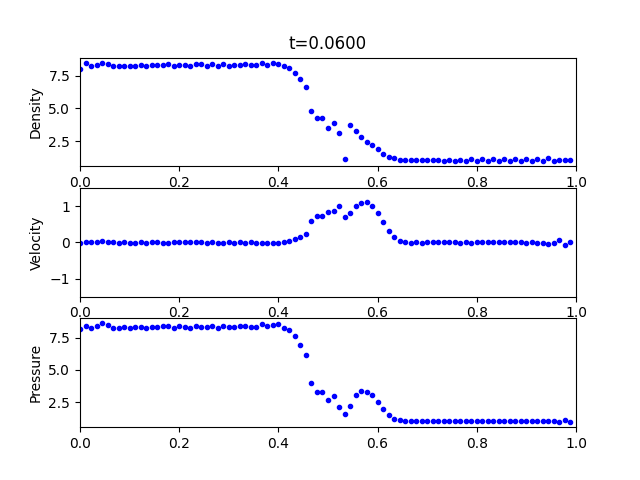}
  \end{minipage}
  }
  \subfigure{
  \begin{minipage}[b]{.4\linewidth}
  \includegraphics[width=1.1\textwidth,height=0.9\textwidth]{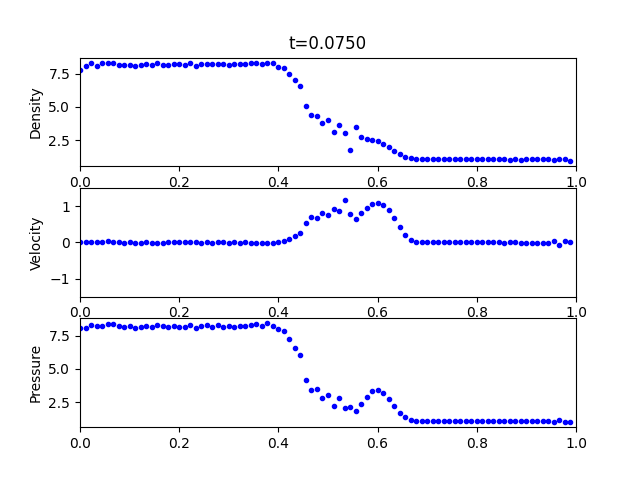}
  \end{minipage}
  }
  \label{fig:fig3}
  \caption{Neural network trained solutions for the Euler's system, by using its integral form. }
\end{figure}

Considering that the macroscopic quantities (the density, velocity and pressure) \LL{carry} quite different physical meanings, we use $3$ independent networks to compute each of them. Each network contains one hidden layer with $1024$ hidden units. The dimensionality of the output vector is now $100$, instead of $300$ used in the above single-model methodology. 
We show the numerical results in Figure \ref{fig:multi} below. Figure \ref{fig:loss} shows the plot of loss value \LL{with respect to} the number of training iteration at different times, \LL{indicating a fast decay of the loss function.}
    
Notice that this is a shock problem, which brings challenges to several general finite difference schemes for hyperbolic systems, unless one uses carefully designed robust solvers \cite{Euler3,VanLeer}. Thus our proposed method has shown its advantages on efficiency and accuracy. 

\begin{figure}[!h]
  \centering
  \subfigure{
  \begin{minipage}[b]{.4\linewidth}
  \includegraphics[width=1.1\textwidth,height=0.9\textwidth]{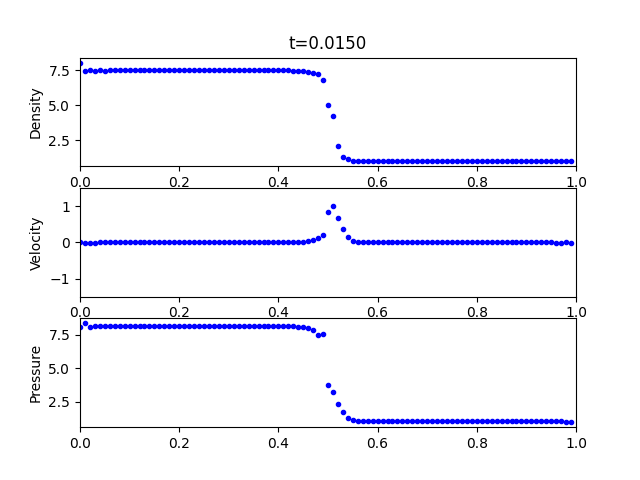}
  \end{minipage}
  }
  \subfigure{
  \begin{minipage}[b]{.4\linewidth}
  \includegraphics[width=1.1\textwidth,height=0.9\textwidth]{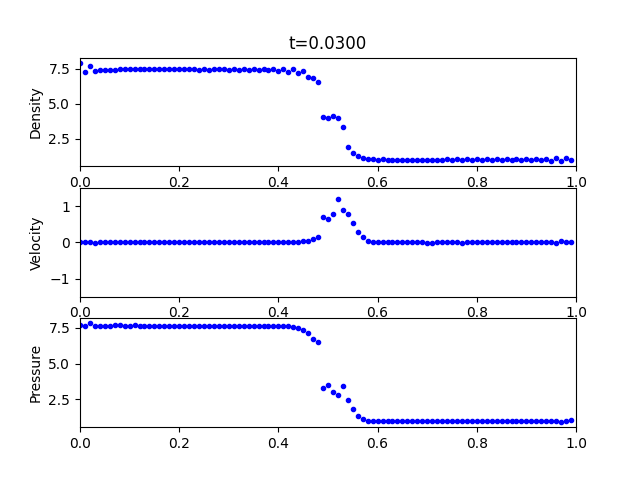}
  \end{minipage}
  }
  \subfigure{
  \begin{minipage}[b]{.4\linewidth}
  \includegraphics[width=1.1\textwidth,height=0.9\textwidth]{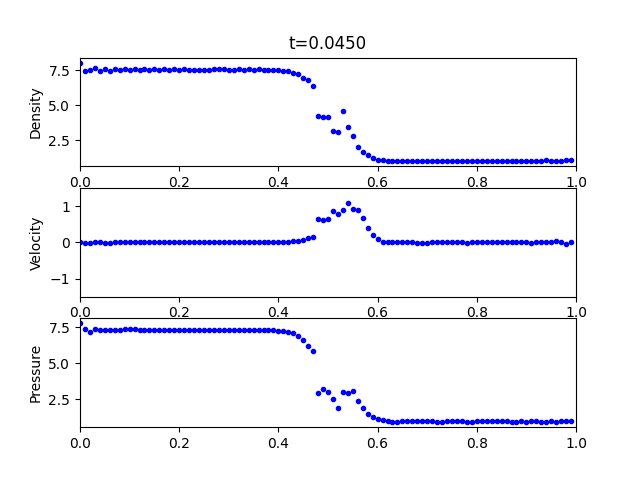}
  \end{minipage}
  }
  \subfigure{
  \begin{minipage}[b]{.4\linewidth}
  \includegraphics[width=1.1\textwidth,height=0.9\textwidth]{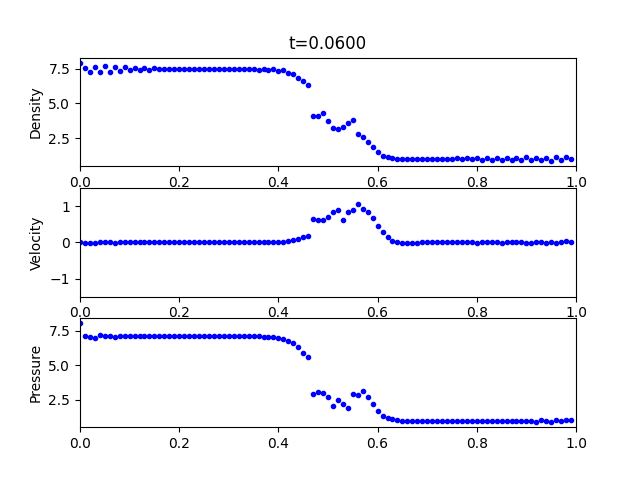}
  \end{minipage}
  }
  \subfigure{
  \begin{minipage}[b]{.4\linewidth}
  \includegraphics[width=1.1\textwidth,height=0.9\textwidth]{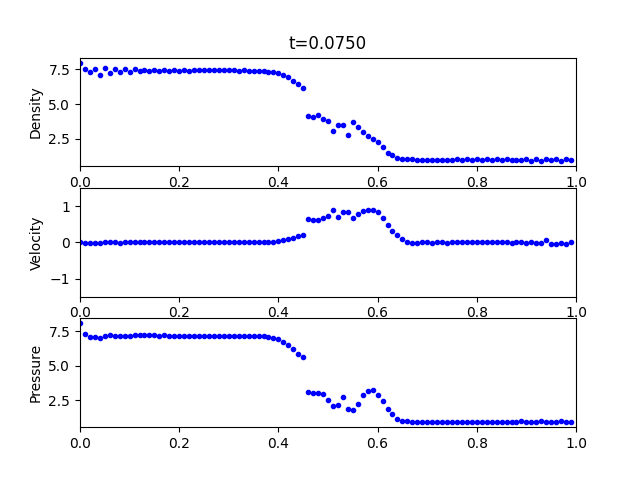}
  \end{minipage}
  }
  \caption{Neural network trained solutions for the Euler's equations with multi-network model, by using its integral form. }
  \label{fig:multi}
\end{figure}

\begin{figure}[!h]
 \label{fig:diff}
   \centering
   \subfigure{
   \begin{minipage}[b]{.4\linewidth}
   \includegraphics[width=1.1\textwidth,height=0.9\textwidth]{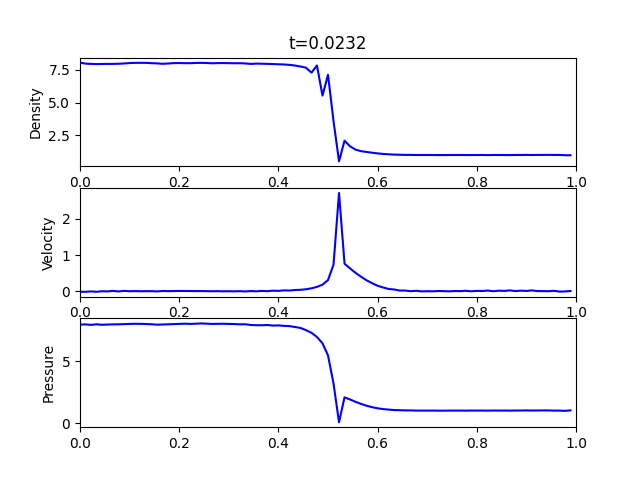}
   \end{minipage}
   }
   \subfigure{
   \begin{minipage}[b]{.4\linewidth}
   \includegraphics[width=1.1\textwidth,height=0.9\textwidth]{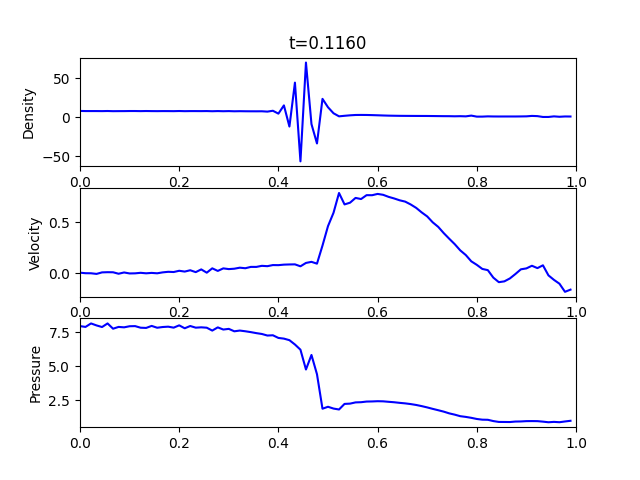}
   \end{minipage}
   }
   \caption{Neural network trained solutions for the Euler's system, by using the differential form.}
 \end{figure}

\begin{figure}[!h]
\centering
\subfigure{
  \begin{minipage}[b]{.3\linewidth}
\includegraphics[width=1.1\textwidth,height=0.9\textwidth]{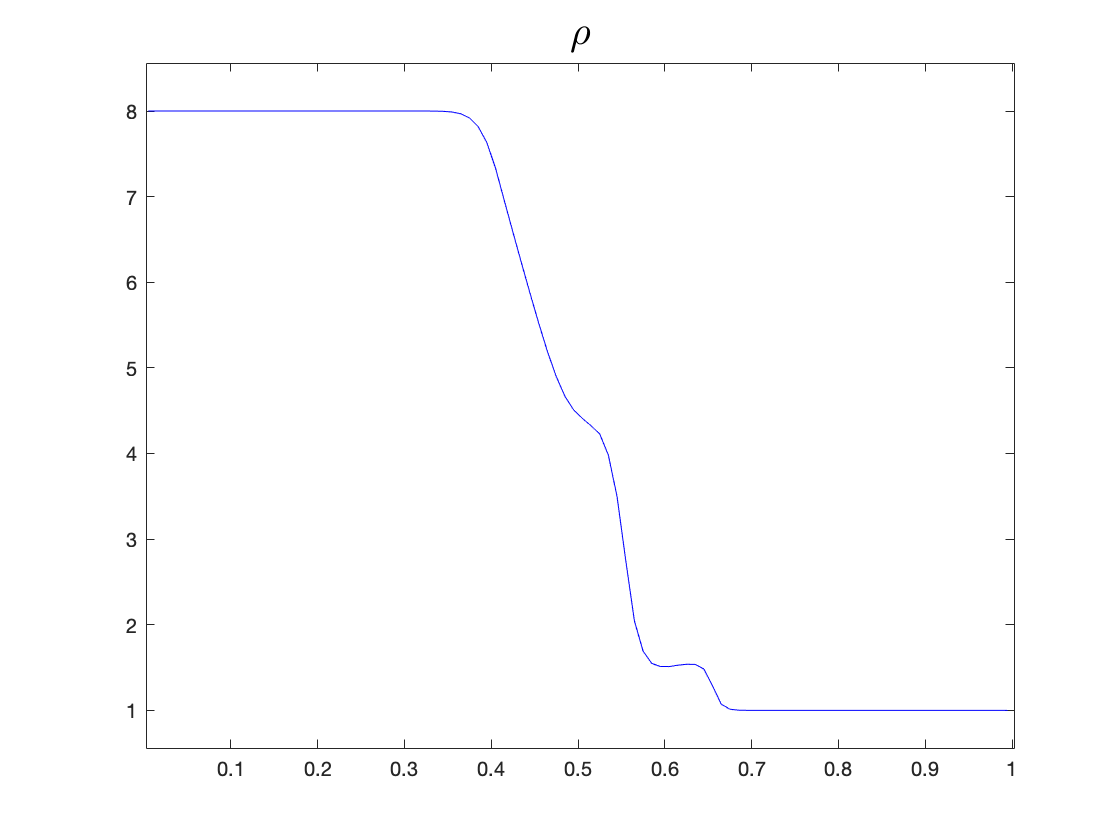}
\end{minipage}
}
\subfigure{
  \begin{minipage}[b]{.3\linewidth}
\includegraphics[width=1.1\textwidth,height=0.9\textwidth]{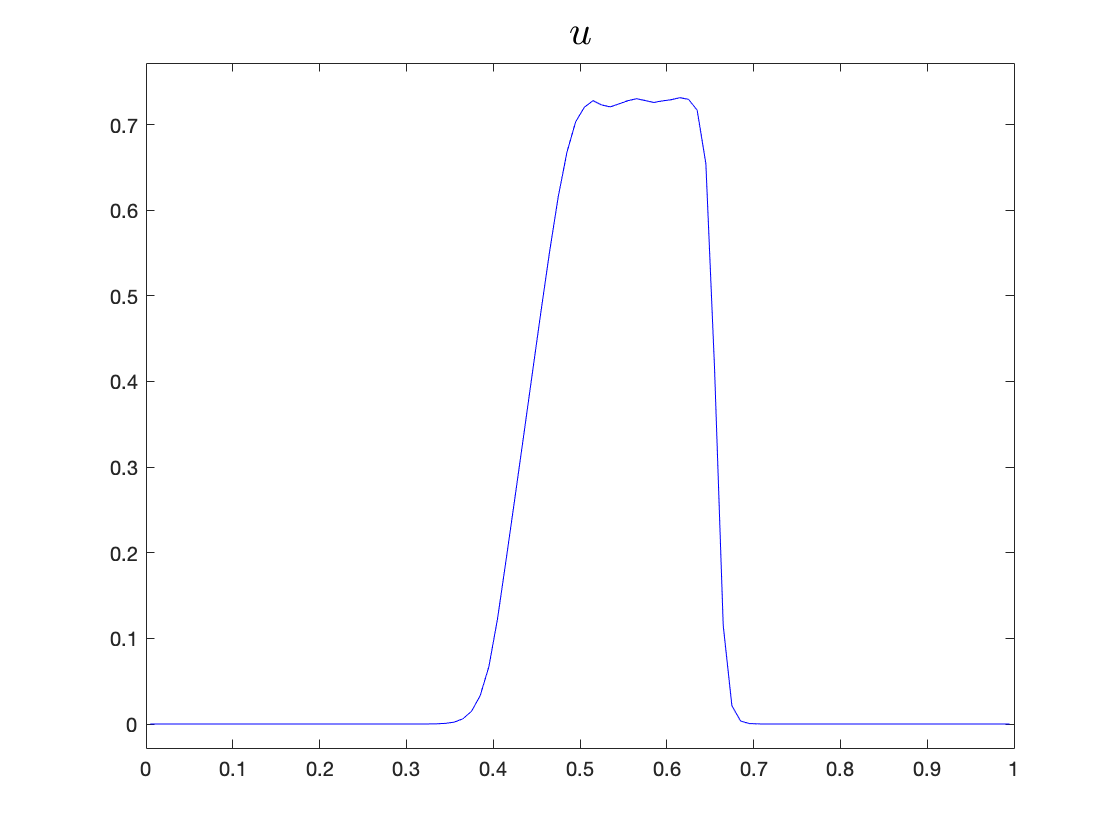}
\end{minipage}
}
\subfigure{
  \begin{minipage}[b]{.3\linewidth}
\includegraphics[width=1.1\textwidth,height=0.9\textwidth]{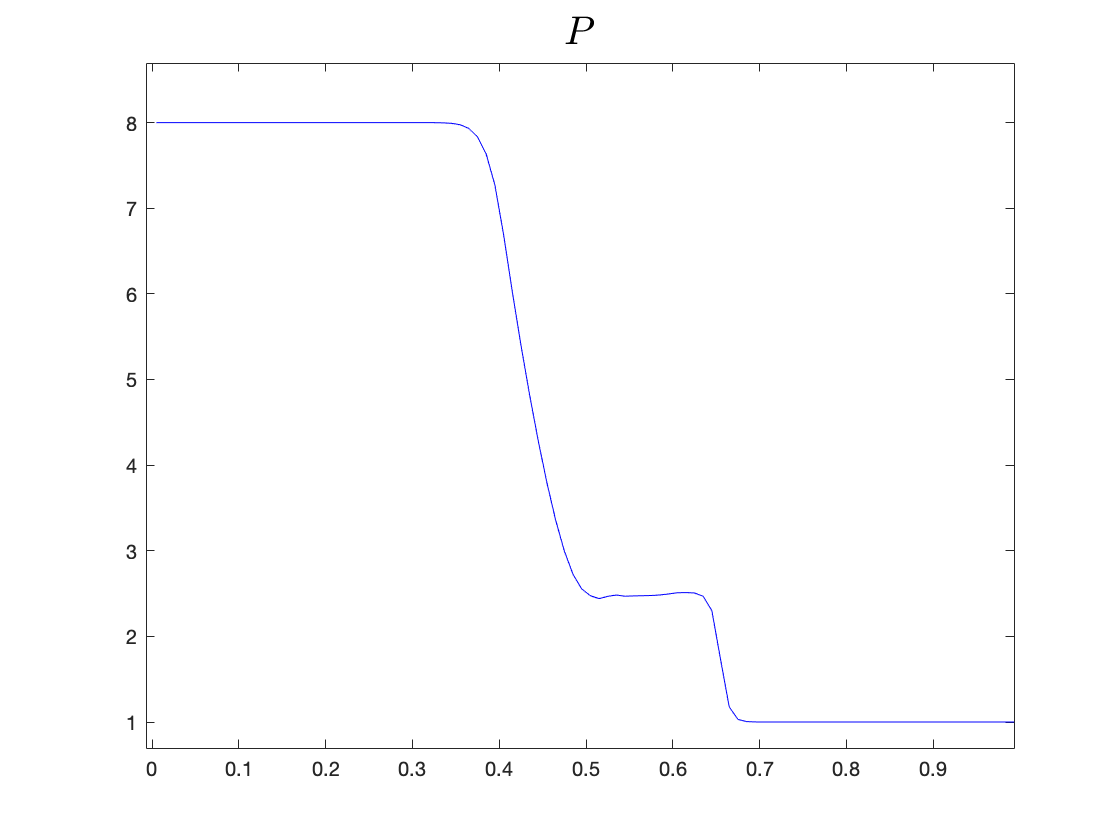}
\end{minipage}
}
\subfigure{
  \begin{minipage}[b]{.3\linewidth}
\includegraphics[width=1.1\textwidth,height=0.9\textwidth]{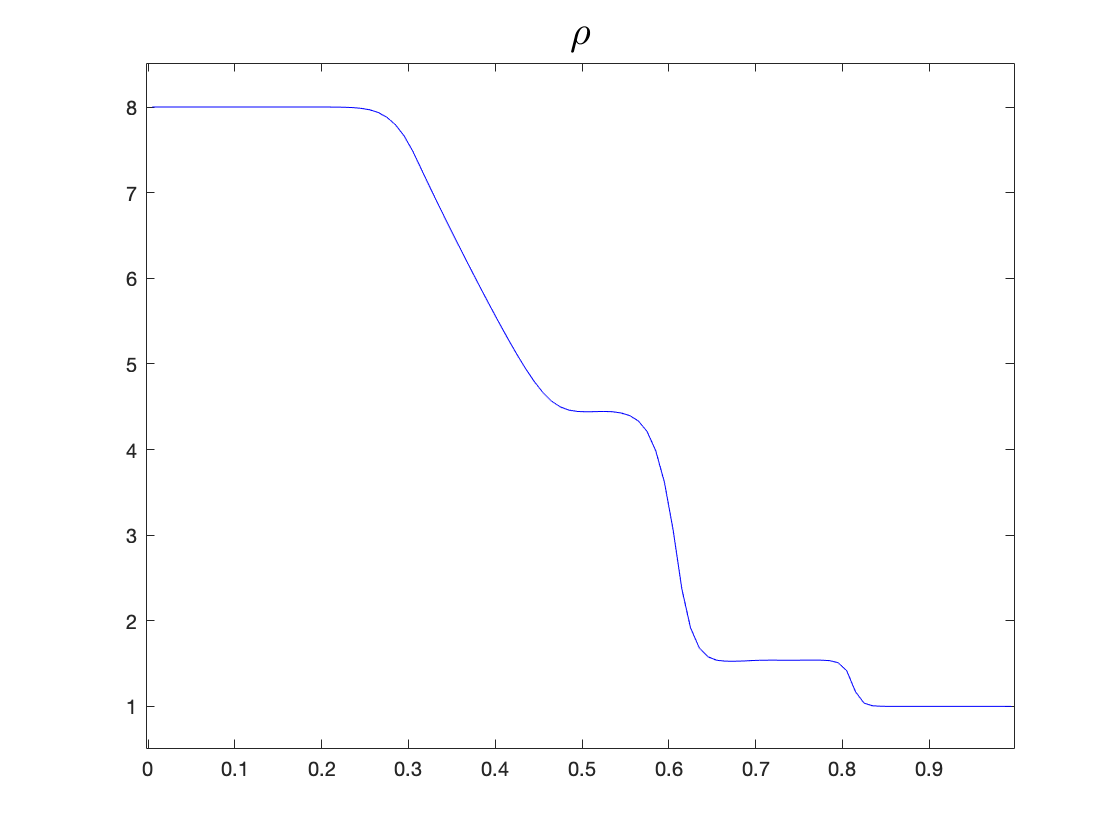}
\end{minipage}
}
\subfigure{
  \begin{minipage}[b]{.3\linewidth}
\includegraphics[width=1.1\textwidth,height=0.9\textwidth]{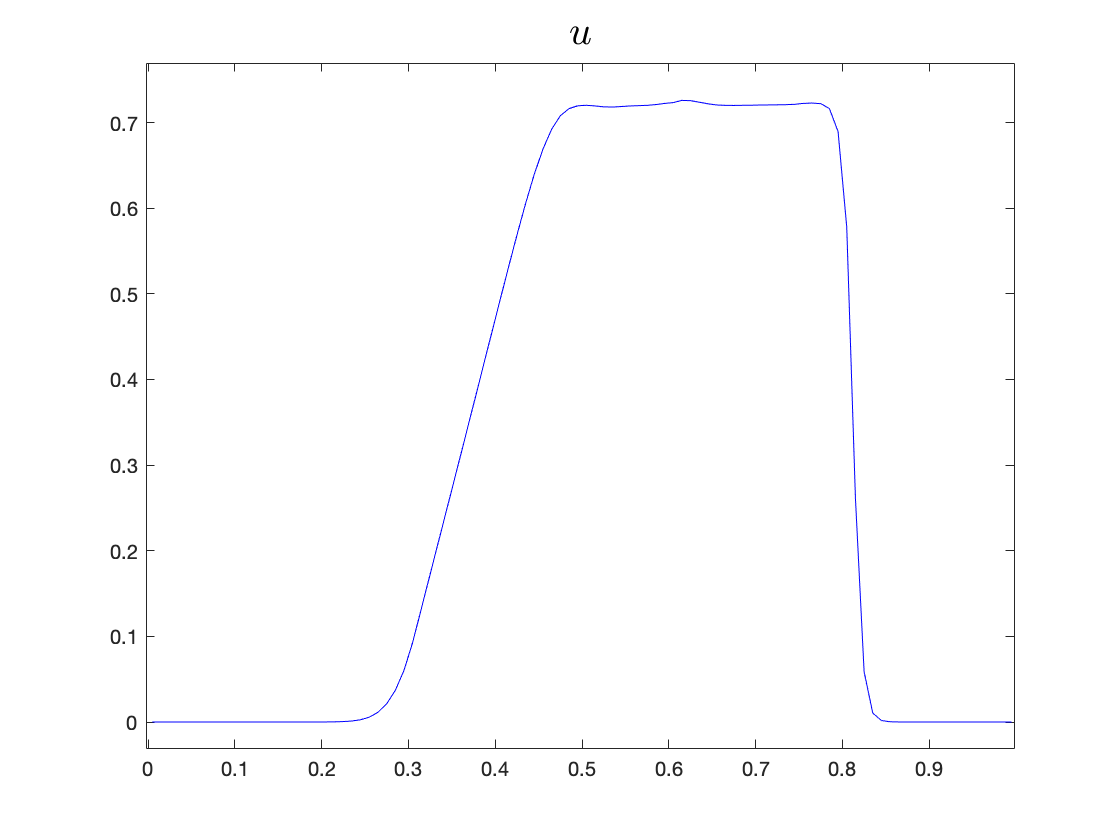}
\end{minipage}
}
\subfigure{
  \begin{minipage}[b]{.3\linewidth}
\includegraphics[width=1.1\textwidth,height=0.9\textwidth]{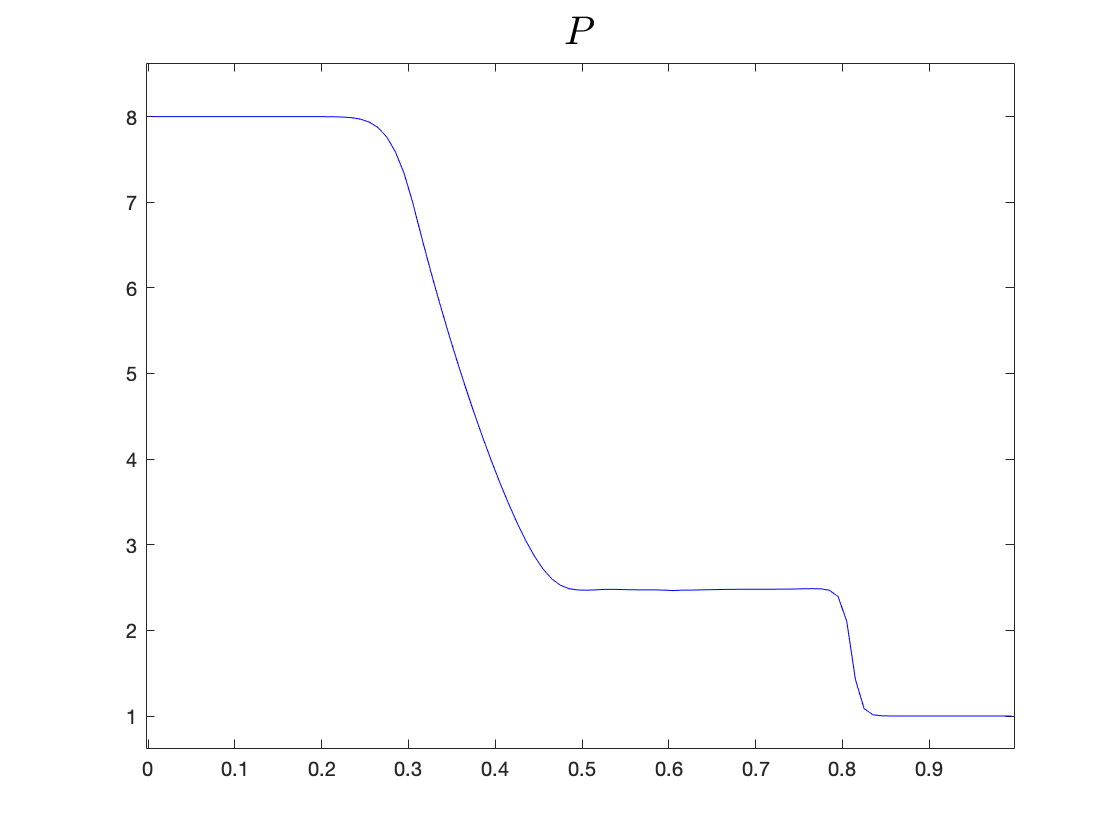}
\end{minipage}
}
\caption{Numerical results obtained by the second-order MUSCL scheme. Output time: $t=0.075$ (above) and $t=0.15$ (bottom).}
\label{fig:ref}
\end{figure}

As a comparison to the integral form of the proposed framework, we conduct experiments by using the original differential form of the Euler's equations. All hyper-parameters and the network structure are the same as used in the integral form. The simulation results are shown in Figure \ref{fig:diff}. 
From above figures, one can observe that the performance of the integral forms are superior than the differential form in both single-network model and multi-network model. This is consistent with our expectations that the approximation used in integral form is of higher-order than the differential form. \LL{In Figure \ref{fig:ref}, as a reference, we show the numerical solution obtained by an accurate second-order MUSCL scheme \cite{VanLeer}, a highly robust solver that has been popularly used for studying the conservation laws.}

\begin{figure}[!h]
\label{fig:loss}
  \centering
  \subfigure{
  \begin{minipage}[b]{.4\linewidth}
  \includegraphics[width=1.1\textwidth,height=0.8\textwidth]{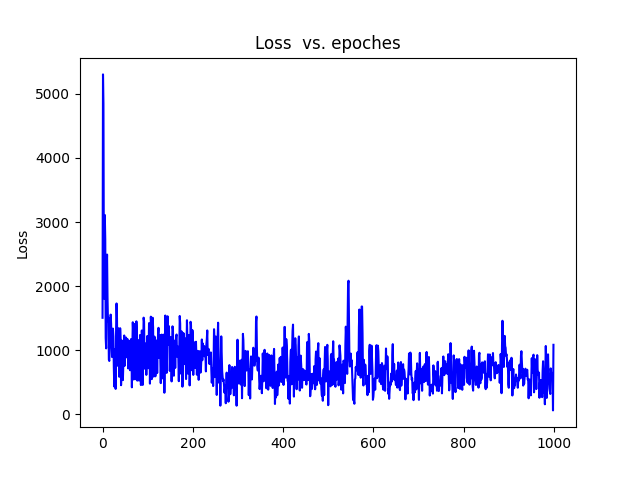}
  \end{minipage}
  }
  \subfigure{
  \begin{minipage}[b]{.4\linewidth}
  \includegraphics[width=1.1\textwidth,height=0.8\textwidth]{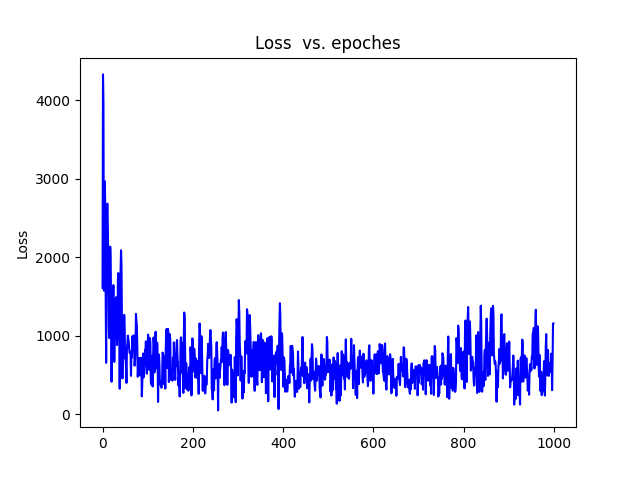}
  \end{minipage}
  }
  \label{fig:fig8}
  \caption{Plot of the loss values with respect to the training epoch.}
\end{figure}

\section{Conclusion and Future Work}
\label{sec:Con}

In this paper, we propose a novel finite difference net, a deep recurrent framework for solving time-dependent evolution equation. It is a new deep learning method for tackling time-dependent PDEs,  motivated by the similarities between finite difference schemes and deep neural networks. It requires no additional training data except the initial or boundary conditions, which is essential to the uniqueness of the solution. By taking integration on both sides of the model, we also use its integral form, 
which is straightforward to compute and more accurate than using the differential form. Numerical experiments are conducted on one-way wave equation and compressible Euler's equations to demonstrate the performance of our finite difference nets.

\LL{For future work, there remains a lot to explore, such as employ our method to solve more general type of and high-dimensional PDEs. How to improve the accuracy of the neural network approximations and develop mathematical theory is another mystery. It is also possible that by using more advanced sampling methods, rather than the uniform sampling for the spatial variable used here, the performance of the finite difference nets may be further improved.}

\bibliographystyle{siam}
\bibliography{FD_Net.bib}

\end{document}